\documentclass[11pt]{article}

\usepackage[latin1]{inputenc}

\usepackage{amsfonts,amsmath,amssymb,amsthm,graphicx,epsfig,float}
\usepackage[T1]{fontenc}

\usepackage[english,francais]{babel}

{
  \newtheorem{theoreme}{Th\'eor\`eme}
  \newtheorem*{theoreme*}{Th\'eor\`eme}
  \newtheorem{lemme}[theoreme]{Lemme}

  \newtheorem{definition}{D\'efinition}
  \newtheorem{proposition}[theoreme]{Proposition}
\newtheorem*{corollaire*}{Corollaire}
\newtheorem*{proposition*}{Proposition}
\theoremstyle{remark}
  \newtheorem*{remarque*}{Remarque}
}

\newcounter{ex}

\newenvironment{rem*}{
  \noindent\textbf{Remarque. }}{}

\newcommand{\Cc}{\mathbb{C}}
\newcommand{\Nn}{\mathbb{N}}

\newcommand{\Pp}{\mathbb{P}}

\newcommand{\Dcal}{\mathcal{D}}

\newcommand{\Hcal}{\mathcal{H}}
\newcommand{\Lcal}{\mathcal{L}}
\newcommand{\II}{\mathcal{I}}
\newcommand{\Qcal}{\mathcal{Q}}
\newcommand{\Scal}{\mathcal{S}}

\title{{\bf Un critère de laminarité locale en dimension quelconque}}
\author{Henry de Thélin}
\date{}

\begin{document}
\maketitle


\def\figurename{{Fig.}}%
\def\proofname{Preuve}
\def\contentsname{Sommaire}%

\begin{abstract}

Nous montrons qu'une suite de sous-ensembles analytiques lisses de
dimension $s$ de la boule unité de $\Cc^l$, dont la courbure est
contrôlée par le volume, converge vers une lamination de dimension $s$
dans un sens faible.

\end{abstract}

\selectlanguage{english}
\begin{center}
{\bf{A criterion of local laminarity for all dimensions}}
\end{center}

\begin{abstract}

We show that a sequence of smooth analytic subsets of dimension $s$ of
the unit ball of $\Cc^l$, for which the curvature is bounded by the
volume, converges to a lamination of dimension $s$ in a weak sense.

\end{abstract}

\selectlanguage{francais}

Mots-clefs: courants, laminarité.

Classification: 32U40, 32H50.

\section*{{\bf Introduction}}
\par

Dans cet article, on s'intéresse à des limites de sous-ensembles
analytiques $M_n$ de dimension $s$ de la boule unité de $\Cc^l$. La
question est de savoir si celles-ci ont conservé un certain caractère
analytique.

Lorsque le volume de $M_n$ est uniformément borné, on sait que c'est
le cas. En effet, quitte à extraire une sous-suite, $M_n$ converge
vers un sous-ensemble analytique de $B$: c'est le théorème de Bishop
(voir \cite{Bi}).

Quand le volume n'est plus majoré et que $s=1$, $l=2$, cette question a été résolue dans
\cite{Det1}. Plus précisément, dans cet article on a montré le:

\begin{theoreme*} Voir \cite{Det1}

Soit $C_n$ une suite de courbes analytiques lisses de la boule unité
$B$ de $\Cc^2$.

On note $A_n$ l'aire de $C_n$, $G_n$ le genre de $C_n$ et on suppose
que $T_n=\frac{[C_n]}{A_n}$ converge vers un $(1,1)$-courant positif
fermé $T$ de $B$ (toujours possible quitte à extraire une sous-suite).

Alors, si $G_n=O(A_n)$, $T$ est laminaire.
\end{theoreme*}

Un courant positif de bidimension $(1,1)$ est laminaire s'il s'écrit localement comme
une intégrale de courants d'intégration sur une famille de disques
disjoints, hors d'un ensemble négligeable (voir \cite{BLS}).

Par ailleurs, les exemples de Wermer (voir \cite{DuvSi}) permettent de construire des
suites de courbes analytiques lisses $C_n$ de la boule unité $B$ dont
le genre croît plus vite que l'aire aussi lentement que l'on veut et pour
lesquelles les limites ne contiennent aucun disque analytique. 

Signalons aussi que le théorème précédent est une version locale de résultats de nature
globale dans $\Pp^2(\Cc)$ (i.e. où les courbes $C_n$ n'ont pas de
bord) qui sont dus à E. Bedford, M. Lyubich et J. Smillie (voir
\cite{BLS}) et R. Dujardin (voir \cite{Du1}). Par ailleurs, tous ces
résultats ont permis de montrer que certains courants issus de la
dynamique holomorphe et méromorphe sont laminaires (voir \cite{BLS},
\cite{Du1} et \cite{Det2}).

$ $

L'objectif de cet article est de traiter le cas où les dimensions $s$
et $l$ sont plus grandes. Pour cela, on remplace tout d'abord la
notion de courants laminaires par celle de courants tissés (voir
\cite{Di}). Un courant positif de bidimension $(s,s)$ est tissé s'il s'écrit localement comme
une intégrale de courants d'intégration sur une famille de boules de
dimension $s$, hors d'un ensemble négligeable (voir le paragraphe
\ref{tissé} pour plus de détails). Ensuite, on remplacera la
notion de genre par celle d'une courbure que l'on va définir maintenant.

On  considère une suite $M_n$ de sous-ensembles analytiques lisses de dimension
$s$ de la boule unité $B$ de $\Cc^l$. Dans toute la suite, on considèrera $\Cc^l$ comme sous-ensemble de
$\Pp^l(\Cc)$ et on notera $\widetilde{G}(l-s,l)$ l'ensemble des plans
complexes de dimension $l-s$ dans $\Pp^l(\Cc)$ (voir par exemple \cite{Ch} p. 165). Les ensembles
$$\widetilde{M_n}=\{ (x, \Dcal) \mbox{, } x \in M_n \mbox{, } x \in \Dcal \in
  \widetilde{G}(l-s,l) \mbox{ et } \Dcal \mbox{ pas
  transverse à } T_x M_n \}$$
vivent dans l'espace produit $\Cc^l \times
  \widetilde{G}(l-s,l)$. C'est par définition le volume de $\widetilde{M_n}$
  que l'on appellera courbure de $M_n$.

Le but de cet article est alors de démontrer le critère suivant:

\begin{theoreme}

Soit $M_n$ une suite de sous-ensembles analytiques lisses de dimension $s$ de
la boule unité $B$ de $\Cc^l$.

On suppose que $T_n=\frac{[M_n]}{\mbox{Volume}(M_n)}$ converge vers un courant positif
fermé $T$ de bidimension $(s,s)$ de $B$ (toujours possible quitte à
extraire une sous-suite).

Alors, si $\mbox{Volume}(\widetilde{M_n})=O(\mbox{Volume}(M_n))$, $T$
est tissé. De plus si $s=l-1$ alors $T$ est laminaire.

\end{theoreme}

Signalons que les exemples précédents de Wermer montrent
que ce théorème est optimal. Par ailleurs, une version globale du
critère (i.e. dans le cas où les $M_n$ n'ont pas de bord) a été démontrée par T.-C. Dinh
dans \cite{Di}. Enfin, de même que le critère de T.-C. Dinh, notre
théorème a pour vocation de démontrer que certains
courants issus de la dynamique holomorphe et méromorphe sont tissés.

$ $

Voici maintenant le plan de ce texte. Après un premier paragraphe
consacré à des préliminaires, le second traitera le cas d'une suite de
courbes dans $\Cc^l$ (i.e. $s=1$ et $l$ quelconque). Enfin, dans le
dernier paragraphe, on démontrera le critère précédent pour toutes les dimensions
de sous-variétés analytiques.

$ $

{{\bf Remerciements:}} Je remercie T.-C. Dinh pour les discussions
fructueuses que nous avons eues au sujet de cet article, ainsi que pour ses encouragements. 

\section{{\bf Préliminaires}}{\label{tissé}}

Dans ce paragraphe, nous allons rappeler les notions de courants tissés
(voir \cite{Di}) et courants laminaires (voir par exemple \cite{BLS}, \cite{C},
\cite{Det1} et \cite{Du1}).

Considérons un ouvert $\Omega$ de $\Cc^l$ et $T$ un courant
positif de bidimension $(s,s)$.

\begin{definition}

Le courant $T$ est uniformément tissé de dimension $s$ dans $\Omega$ si pour tout
$x \in \mbox{Supp}(T) \cap \Omega$, il existe un polydisque $B$, un
ouvert $U$ de $B$ contenant $x$ et une constante $c(B,T)$ tels que:
$$T_{|U}= \int_{\Gamma \in \mathcal{G}} [\Gamma \cap U] d \lambda
(\Gamma).$$
Ici $ \mathcal{G}$ est l'ensemble des sous-ensembles analytiques
irréductibles de
dimension $s$ de $B$ de masse inférieure à $c(B,T)$ et $\lambda$ est une mesure sur cet
espace compact.

\end{definition}

Remarquons qu'a priori les $\Gamma$ ne sont pas supposés
disjoints. Cependant, si pour tous $\Gamma$ et $\Gamma'$ du support
de $\lambda$ on a $\Gamma \cap \Gamma' = \varnothing$ ou
$\mbox{Dimension}( \Gamma \cap \Gamma')=s$, on dira que $T$
est uniformément laminaire de dimension $s$.  

De façon analogue aux courants laminaires, on peut maintenant définir:

\begin{definition}

Un courant $T$ est tissé (respectivement laminaire) de dimension $s$ dans $\Omega$ s'il existe une suite
d'ouverts $\Omega_i \subset \Omega$ avec $||T||(\partial \Omega_i)=0$
et une suite croissante $(T_i)_{i \geq 0}$, $\mbox{ }T_i$ uniformément
tissés (respectivement uniformément laminaires) de dimension $s$ dans $\Omega_i$ tels que $\lim_{i \rightarrow \infty}
T_i=T$.

\end{definition}

\section{{\bf Le cas d'une suite de courbes dans $\Cc^l$}}{\label{dim1}}

Considérons une suite $C_n$ de courbes analytiques lisses de la boule unité
de $\Cc^l$ telle que le volume de $\widetilde{C_n} \subset \Cc^l
\times \widetilde{G}(l-1,l)$ soit borné par $O( \mbox{Volume}(C_n))$.

L'objectif de ce paragraphe est de montrer que les valeurs d'adhérence de
$\frac{[C_n]}{\mbox{Volume}(C_n)}$ sont tissées (ou laminaires si
$l=2$).

Voici le plan de la démonstration. Dans un premier paragraphe, nous
allons montrer que le contrôle sur le volume de $\widetilde{C_n}$
implique que
$$ \int_{C_n}KdV \geq - O( \mbox{Volume}(C_n)),$$
où $K$ est la courbure de Gauss de $C_n$. On utilisera pour cela un
résultat de R. Langevin et T. Shifrin (voir \cite{LS}). Dans le second
paragraphe, nous verrons que ce contrôle de la courbure implique que
le genre de $C_n$ est en $O(\mbox{Volume}(C_n))$ quitte à réduire un
peu la boule $B$. Enfin, en utilisant \cite{Det1}, nous en déduirons
que les limites de $\frac{[C_n]}{\mbox{Volume}(C_n)}$ sont tissées (ou
laminaires si $l=2$).

\subsection{{\bf Contrôle de la courbure de Gauss des courbes $C_n$}}{\label{2.1}}

Tout d'abord d'après le théorème $4.3$ de \cite{LS} on sait que
$$-\frac{1}{\pi}  \int_{C_n} K dV= \int_{G(l-1,l)} n(C_n,H) dH,$$
où $n(C_n,H)= \# \{ z \in C_n \mbox{, } T_z C_n \subset H \}$ (compté
avec multiplicité) et
$G(l-1,l)$ est l'ensemble des hyperplans de $\Cc ^l$.

Maintenant, si on fixe un hyperplan $H$ de $\Cc^l$ (par exemple
$z_1=0$) alors les hyperplans affines $z_1=a$ (avec $a \in \Cc$)
décrivent une droite $L$ dans $\widetilde{G}(l-1,l)$. De plus $n(C_n,H)$ est le nombre d'intersection de $L$ avec
la projection  de $\widetilde{C_n}$ sur
$\widetilde{G}(l-1,l)$. Autrement dit $\int_{G(l-1,l)} n(C_n,H) dH$
est majoré par $ \mbox{Volume}(\widetilde{C_n})$.

En combinant les relations obtenues, on a donc bien
$$ \int_{C_n} K dV \geq - O( \mbox{Volume}(C_n)).$$

\subsection{{\bf Passage du contrôle de la courbure à celui du genre}}

Dans ce paragraphe, nous allons montrer la

\begin{proposition}{\label{prop2}}

Soit $C$ une courbe lisse de la boule unité $B$ de $\Cc^l$. On suppose
que $C$ n'a pas de bord dans la boule $B$ et on munit $C$ de la métrique hermitienne induite. Alors, on a:

$$\mbox{Genre}(C \cap \rho B) \leq c( \rho) \left( \mbox{Volume}(C) -
 \int_C KdV \right)$$
pour tout $0 < \rho < 1$.

\end{proposition}

\begin{proof}

Considérons trois boules concentriques $\rho B$, $\rho'B$ et $B$
(où $\rho'$ est générique et strictement compris entre $\rho$ et $1$).

L'idée va être de transformer les composantes de bord de $C \cap
\rho'B$ en des courbes géodésiques par morceaux. On obtiendra ainsi
une surface $\widetilde{C}$ à bord géodésique par morceaux. Le calcul
du genre de $\widetilde{C}$ passera alors par une minoration de
$\chi(\widetilde{C})$ qui sera obtenue par la formule de Gauss-Bonnet.

\subparagraph{{\bf 1) Transformation du bord de $C$}}

$ $

Soit $b$ une composante de bord de $C \cap \rho'B$. Si la longueur de $b$ est
inférieure à $\epsilon$ ($\epsilon << \min(1 - \rho',
\rho'-\rho)$), on note $\widetilde{b}$ la plus petite courbe dans
la classe d'homotopie de $b$ (dans $\overline{C}$). Cette courbe $\widetilde{b}$ est proche
de $b$. En effet, en terme d'homologie $\widetilde{b}-b$ est le bord
de $S$. Mais si $\omega$ désigne la forme kählérienne de $\Cc^l$, on a
$\omega= d \lambda$ dans $B$ d'où:
$$\mbox{Aire}(S) = \int_S \omega = \int_S d \lambda = \int_{\partial
  S} \lambda \leq \| \lambda \| (
  \mbox{Longueur}(\widetilde{b})+\mbox{Longueur}(b)) \leq 2 \epsilon \| \lambda \|  ,$$
car la longueur de $\widetilde{b}$ est inférieure à $\epsilon$.

La courbe $\widetilde{b}$ reste donc dans un $ \epsilon^{1/3}$-voisinage de
$b$ par le théorème de Lelong (voir \cite{L}). En particulier
$\widetilde{b}$ est une géodésique (elle ne touche pas le bord de
$B$). De plus, en utilisant de nouveau le théorème de Lelong, on
constate que $S$ reste dans un $2 \epsilon ^{1/3}$-voisinage de $\partial
\rho'B$.

Traitons maintenant le cas d'une composante de bord $b$ de $C \cap
\rho'B$ de longueur au moins $\epsilon$. 

On découpe $b$ en $\left[ \frac{L}{\epsilon} \right] + 1$ morceaux de
longueur inférieure ou égale à $\epsilon$. Fixons une de ces
composantes connexes. Dans la classe d'homotopie de cette
composante avec les deux extrémités fixées on choisit
la plus petite courbe. Comme elle est de longueur inférieure à
$\epsilon$, elle reste dans un $\epsilon$-voisinage de $\partial
\rho' B$ et elle est donc une géodésique. En recommençant le
procédé avec toutes les composantes connexes de $b$, on obtient des
géodésiques qui une fois réunies forment une courbe $\widetilde{b}$
homotope à $b$ et géodésique par morceaux. Remarquons que
$\widetilde{b}$ vit dans un $\epsilon$-voisinage de $b$ et que par les
mêmes arguments que précédemment, si $\widetilde{b} -b=\partial S$,
$S$ reste dans un $ \epsilon^{1/3}$-voisinage de $b$ (utiliser les
arguments pour chaque composante connexe du découpage de $b$).

Dans la suite, on notera $\Scal$ l'ensemble des sommets qui sont aux
extémités des composantes connexes qui découpaient les bords de
longueur supérieure à $\epsilon$.

Faisons un bilan de ce que l'on a fait: à toute composante de bord $b$
de $C \cap \rho'B$, on a associé une géodésique $\widetilde{b}$ (si
la longueur de $b$ était inférieure à $\epsilon$), ou une géodésique
par morceaux $\widetilde{b}$ (si la longueur de $b$ était supérieure à
$\epsilon$). Le bord de $C \cap \rho'B$ est homologue à cette union
de géodésiques et géodésiques par morceaux $\cup \widetilde{b}$. Enfin $\cup \widetilde{b} - \partial(C \cap \rho'B)$ est le bord de
quelque chose qui vit dans un $2 \epsilon^{1/3}$-voisinage de $\partial
\rho'B$.

Maintenant, si on découpe $C$ suivant les courbes $\widetilde{b}$, on
obtient un certain nombre de composantes connexes. On notera
$\widetilde{C}$ l'union de ces composantes connexes qui rencontrent $C \cap \rho B$.

\subparagraph{{\bf 2) Minoration de  $\chi(\widetilde{C})$}}

$ $

Remarquons tout d'abord que le bord de $\widetilde{C}$ est constitué
uniquement de morceaux des $\widetilde{b}$ construits précédemment. En
effet, si cela n'était pas le cas, on pourrait construire une courbe
$\gamma$ dans la surface $C$ qui joindrait un point de $\rho B$ à un
point de $\partial B$ sans jamais toucher un bord $\widetilde{b}$
(toute composante connexe de $\widetilde{C}$ rentre dans $\rho B$). En
terme d'intersection, on aurait donc $\gamma. (\cup
\widetilde{b})=0$. Par ailleurs $\gamma. \partial(C \cap \rho' B)=
\pm 1$ (le nombre d'intersection produit par le fait que $\gamma$ sort
de $\rho' B$ est l'opposé de celui produit par le fait que $\gamma$
y entre). Pour obtenir une contradiction, il suffit de démontrer le:

\begin{lemme}

Soit $\gamma$ un chemin de $C$.

Si $a$ et $b$ sont homologues dans $C$ ($a-b= \partial S$), et si les
extrémités de $\gamma$ ne sont pas dans $S$, on a $\gamma.a=\gamma.b$.

\end{lemme}

\begin{proof}

Si $\gamma$ rentre dans un simplexe de $S$, il doit en sortir. Mais le
nombre d'intersection produit par le fait d'entrer dans un simplexe
est l'opposé de celui produit par le fait d'en sortir. D'où $\gamma.(a-b)=0$.

\end{proof}

Le bord de $\widetilde{C}$ est donc constitué uniquement de morceaux
géodésiques inclus dans $\cup \widetilde{b}$. En particulier il est
inclus dans un $2 \epsilon^{1/3}$-voisinage de $\partial \rho' B$.

Pour arriver à estimer $\chi(\widetilde{C})$, on doit faire une autre
remarque sur les courbes $\widetilde{b}$: si on se place sur un
morceau lisse connexe $m$ de $\cup \widetilde{b}$ privé des points
multiples, on va voir que $\widetilde{C}$ se trouve au plus d'un côté
de $m$. En effet dans le cas contraire, on pourrait construire deux
courbes $\gamma_1$, $\gamma_2$ sur $\widetilde{C}$, qui partent de
$\rho B$ et qui arrivent chacune d'un côté de $m$. En pertubant un peu
la situation, on peut donc produire une courbe $\gamma$ de $C$ qui
joint deux points de $\rho B$ telle que $\gamma. \cup \widetilde{b}=
\pm 1$. Mais $\gamma. \partial (C \cap \rho'B)=0$ (car on sort de
$\rho' B$ autant de fois qu'on y rentre), ce qui contredit le lemme
précédent.

$ $

Passons maintenant à la minoration de $\chi(\widetilde{C})$.

La formule de Gauss-Bonnet nous donne:

$$ \int_{\widetilde{C}}KdV + \int_{\partial \widetilde{C}} k_g=2 \pi
\chi(\widetilde{C})-2 \pi \chi(\partial\widetilde{C})+
\sum (\beta(s) - \pi)$$
où $k_g$ désigne la courbure géodésique. La dernière somme est prise sur les sommets du bord de $\widetilde{C}$ et
$\beta(s)$ désigne l'angle intérieur au sommet $s$.

Remarquons que l'orientation de $\partial \widetilde{C}$ n'est plus
nécessairement la même que celle induite par l'homologie.

Pour minorer $\chi(\widetilde{C})$, il faut étudier les sommets
$s$ du bord de $\widetilde{C}$.

$ $

\underline{1\up{er} Cas}: $s \in \Scal$ (i.e. le sommet est une
extrémité d'une composante connexe qui découpait un bord long).

Si la valence de $s$ (i.e. le nombre d'arêtes qui
arrivent en $s$) est égale à $2$, alors l'angle intérieur peut être
supérieur à $\pi$. Cependant il y a au plus
$\frac{2}{\epsilon}\mbox{Longueur}(\partial (C \cap \rho' B))$ tels
  sommets.

Si la valence est supérieure à $3$ alors elle est au moins $4$
et on modifie un peu la situation comme dans le 2\up{ème} cas.

$ $

\underline{2\up{ème} Cas}: $s \notin \Scal$ et $s$ est à
l'intersection de plusieurs morceaux de géodésiques (intersection
nécessairement transverse).

Dans ce cas, on a un certain nombre de secteurs angulaires et par la
remarque précédente, la surface $\widetilde{C}$ ne peut pas se trouver
dans deux secteurs adjacents. Maintenant, si on considère un secteur
angulaire où $\widetilde{C}$ se trouve, on le modifie comme dans la
figure \ref{Figure2}.

\begin{figure}[c,h]
\begin{center}

\setlength{\unitlength}{1cm}
\begin{picture}(16,3)
\includegraphics{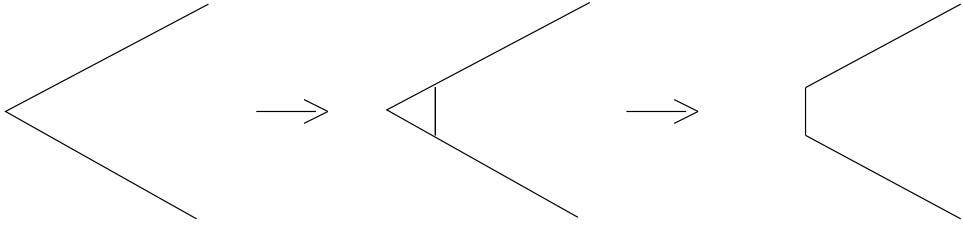}
\end{picture}

\end{center}
\caption{Transformation du bord pour un secteur angulaire{\label{Figure2}}.}
\end{figure}

Cette opération ajoute des sommets mais pour chacun d'entre eux l'angle
intérieur est inférieur ou égal à $\pi$. Grâce à ces transformations
du bord de $\widetilde{C}$, on obtient une nouvelle surface
$\widetilde{C}$ dont le bord est plongé et géodésique par morceaux. Par ailleurs l'angle
intérieur en un sommet du bord de $\widetilde{C}$ est toujours
inférieur à $\pi$ sauf pour au plus $\frac{2}{\epsilon}
\mbox{Longueur}(\partial(C \cap \rho' B))$ d'entre eux. La formule
de Gauss-Bonnet nous donne donc:
$$  \int_{\widetilde{C}} KdV =2
\pi \chi(\widetilde{C}) + \sum(\beta(s)-\pi)$$
où la dernière somme est prise sur l'ensemble des sommets $s$ du bord de
$\widetilde{C}$.

On a donc:
$$ \int_{\widetilde{C}} KdV \leq 2 \pi \chi(\widetilde{C}) +
\frac{2 \pi}{\epsilon}
\mbox{Longueur}(\partial(C \cap  \rho' B)),$$
car $\beta(s) \leq 2 \pi$ si $s \in \Scal$.
Mais d'une part
$$ \int_{\widetilde{C}} KdV \geq  \int_{C} KdV$$
car $K$ est négative et d'autre part 
$$\mbox{Longueur}(\partial(C \cap  \rho' B)) \leq c(\rho)
\mbox{Volume}(C)$$
si $\rho'$ est générique par la formule de la coaire (voir par
exemple \cite{Fe}).

Autrement dit:
$$ \chi(\widetilde{C}) \geq  \frac{1}{2 \pi}  \int_{C} KdV -
c(\rho) \mbox{Volume}(C).$$
Pour finir la démonstration, il reste à contrôler le genre de $C \cap
\rho B$ grâce à la minoration de $\chi(\widetilde{C})$.

On a:
$$\mbox{Genre}(C \cap \rho B) \leq \mbox{Genre}(\widetilde{C})=
\mbox{Nombre de composantes connexes de } \widetilde{C} -
\frac{\chi(\widetilde{C})+p}{2}$$
où $p$ est le nombre de composantes de bord de $\widetilde{C}$.

Comme les composantes connexes de $\widetilde{C}$ ont leur bord dans
un $2 \epsilon^{1/3}$-voisinage de $\partial \rho'B$ et qu'elles entrent
dans $\rho B$ par définition, elles sont en nombre au plus égal à
$c(\rho) \mbox{Volume}(C)$ en utilisant le théorème de
Lelong. Finalement, on a:
$$\mbox{Genre}(C \cap \rho B) \leq c(\rho) \left( \mbox{Volume}(C) - \int_{C} KdV \right),$$
qui est l'inégalité que l'on voulait démontrer.

\end {proof}

\subsection{{\bf Caractère tissé de la limite}}

Quitte à extraire une sous-suite $\frac{[C_n]}{\mbox{Volume}(C_n)}$
converge vers un courant $T$. Le but de ce paragraphe est de démontrer
que $T$ est tissé (ou laminaire si $l=2$). Ces notions étant locales,
il suffit de le démontrer dans $\rho B$ (avec $\rho$ générique compris
entre $0$ et $1$).

Si on combine ce que l'on a obtenu dans les deux paragraphes précédents on a que le genre de $C_n$ dans $\rho B$ est majoré
par $O(\mbox{Volume}(C_n))$. L'idée maintenant va être d'utiliser la
démonstration de \cite{Det1}, afin de prouver le caractère tissé de
$T$.

On a deux possibilités. Soit $T$ est nul dans
$\rho B$ (et le théorème est démontré), soit il existe $D$ une
direction pour laquelle $\pi_{*}(T_{| \rho B}) \neq 0$ (où
$\pi$ est la projection orthogonale sur la droite $D$). Dans toute la suite on
se placera dans ce dernier cas. Si on quadrille le carré $C \subset
D$, centré en $0$, de côté $2$ en $4k^2$ carrés égaux, alors dans
\cite{Det1}, on a démontré que le nombre de bonnes îles dans $C_n \cap
\rho B$
(i.e. graphes au-dessus des carrés du quadrillage) est minoré par $4k^2(1 - \epsilon_k)S_n$. Ici $S_n$
est essentiellement le recouvrement moyen de $C_n \cap \rho B$
au-dessus de $C$ (i.e. $S_n= \frac{1}{\mbox{Volume}(C)} \int_{C_n \cap
\rho B}
\pi^{*} \omega$ où $\omega$ est la forme kählérienne de $\Cc$). Par
ailleurs $\epsilon_k$ est une suite qui tend lentement vers $0$ (et
dans ce texte toute suite de ce type sera notée $\epsilon_k$).

Montrons maintenant que $T$ est tissé grâce à la minoration du nombre
de bonnes îles.

Soit $T_{k,n}$ le courant défini par $T_{k,n}=\frac{1}{\mbox{Volume}(C_n)}
\displaystyle\sum_{ \mbox{ bonnes îles}}[\Gamma]$. Le courant
$T_{k,n}$ peut aussi s'écrire $T_{k,n}=\int [\Gamma] d \nu_{k,n}(\Gamma)$ où
$\nu_{k,n}$ est une mesure sur l'espace métrique compact des graphes
au-dessus des carrés du quadrillage.

Si on note $T_n=\frac{[C_n \cap \rho B]}{\mbox{Volume}(C_n)}$ alors on a:
$$\int{T_{k,n} \wedge \pi^*\omega} \geq (1-\epsilon_k)\int{T_{n} \wedge \pi^*\omega},$$
d'où,
$$\int{(T_{n}-T_{k,n}) \wedge \pi^*\omega} \leq \epsilon_k.$$
La suite de mesures $\nu_{k,n}$ converge vers une mesure $\nu_k$
(quitte à extraire une sous-suite) ce qui implique que $T_{k,n}$
converge vers $T_k= \int [\Gamma] d \nu_k(\Gamma) $ qui est donc uniformément
tissé au-dessus de chaque carré du quadrillage (et uniformément
laminaire si $l=2$). Par ailleurs, on a toujours l'estimée:
$$\int{(T_{| \rho B}-T_{k}) \wedge \pi^*\omega} \leq \epsilon_k,$$
avec $T_{| \rho B}-T_k \geq 0$ par construction.
Si on raffine de plus en plus le quadrillage (i.e. si $k$ augmente),
$T_k$ croît vers un courant $T_{\infty}$ qui est tissé (ou laminaire
si $l=2$). De plus $T_{\infty} \leq T_{| \rho B}$ et $\int{(T_{| \rho B}-T_{\infty}) \wedge \pi^*\omega} \leq
0$.

Maintenant, si on prend une autre direction $D'$ générique par
rapport à $T_{\infty}$ et telle que $\pi_{*}'(T_{| \rho B}) \neq 0$ (où
$\pi'$ désigne la projection associée à $D'$), on construit de
même un courant $T_{\infty}' \leq T_{| \rho B}$ qui est supérieur à
$T_{\infty}$ et qui vérifie  $\int{(T_{| \rho B} -T_{\infty}') \wedge
  \pi'^{*} \omega} = 0$. En itérant ce procédé on finit par avoir $T_{|
  \rho B}=T_{\infty}'$, c'est-à-dire que $T$ est tissé dans $\rho
B$ (ou laminaire si $l=2$).

\section{{\bf Le cas général}}

On considère ici une suite $M_n$ de sous-ensembles analytiques lisses de
dimension $s$ de la boule unité $B$ de $\Cc^l$. Par hypothèse le
volume de 
$$\widetilde{M_n}=\{ (x, \Dcal) \mbox{, } x \in M_n \mbox{, } x \in \Dcal \in
  \widetilde{G}(l-s,l) \mbox{ et } \Dcal \mbox{ pas
  transverse à } T_x M_n \}$$
dans $\Cc^l \times \widetilde{G}(l-s,l)$ est
contrôlé par $O(\mbox{Volume}(M_n))$. Quitte à extraire une sous-suite
$\frac{[M_n]}{\mbox{Volume}(M_n)}$ converge vers un courant $T$ et
nous voulons démontrer que $T$ est tissé (ou laminaire si $s=l-1$).

Voici l'idée de la preuve. Dans un premier temps nous allons trancher $M_n$ avec des plans complexes de dimension $l-s+1$. On obtiendra
ainsi des courbes dans $\Cc^l$. Nous verrons que la majoration
du volume de $\widetilde{M_n}$ impliquera un bon contrôle de la
courbure des courbes obtenues. En particulier, nous pourrons donc
construire beaucoup de disques sur ces courbes en utilisant le
paragraphe précédent. Pour conclure il nous restera alors à utiliser
le théorème de N. Sibony et P. M. Wong (\cite{SW}) pour passer des
disques à des boules de dimension $s$.

\subsection{{\bf Tranchage de $M_n$ par des plans complexes de
    dimension $l-s+1$}}{\label{3.1}}

Dans ce paragraphe, nous allons estimer la courbure des courbes $M_n
\cap \Hcal$ avec $\Hcal \in \widetilde{G}(l-s+1,l)$. On utilisera à
plusieurs reprises la formule de la coaire que l'on peut trouver dans
\cite{Fe} page 258.

Commençons par considérer l'ensemble
$$\widetilde{\widetilde{M_n}}= \{ (x, \Dcal, \Hcal) \mbox{, } x \in
M_n \mbox{, } x \in \Dcal \in \widetilde{G}(l-s,l) \mbox{ et } \Dcal
\mbox{ pas transverse
  à } T_x M_n \mbox{, } \Dcal \subset \Hcal \in \widetilde{G}(l-s+1,l)  \}.$$
Le volume de $\widetilde{M_n}$ étant contrôlé par
$O(\mbox{Volume}(M_n))$ celui de $\widetilde{\widetilde{M_n}}$ l'est
aussi.

Maintenant si on note $\pi_3$ la projection de $\Cc^l \times
\widetilde{G}(l-s,l) \times \widetilde{G}(l-s+1,l)$ sur
$\widetilde{G}(l-s+1,l)$ et que $d \Hcal$ désigne la mesure volume
sur $\widetilde{G}(l-s+1,l)$, on a:
$$\int \mbox{Volume}(\widetilde{\widetilde{M_n}} \cap \pi_3^{-1}(\Hcal)) d \Hcal =
O(\mbox{Volume}(M_n))$$
grâce à la formule de la coaire.

Autrement dit,
$$ \int \mbox{Volume}((M_n)_{\Hcal}) d \Hcal = O(\mbox{Volume}(M_n))$$
où $(M_n)_{\Hcal} \subset \Cc^l \times \widetilde{G}(l-s,l)$ est
défini par:
$$(M_n)_{\Hcal}= \{ (x, \Dcal) \mbox{, } x \in M_n \mbox{, } x \in \Dcal \in
\widetilde{G}(l-s,l) \mbox{, } \Dcal \mbox{ pas transverse à } T_x M_n \mbox{ et }
\Dcal \subset \Hcal \}.$$
En particulier pour $\Hcal$ générique $M_n \cap \Hcal$ est une
courbe lisse et:
$$(M_n)_{\Hcal}= \{ (x, \Dcal) \mbox{, } x \in M_n \cap \Hcal \mbox{, } x \in \Dcal \in
\widetilde{G}(l-s,l) \mbox{, } T_x (M_n \cap \Hcal) \subset \Dcal
\subset \Hcal \}.$$
Le volume de $(M_n)_{\Hcal}$ est donc égal à celui de $\widetilde{(M_n \cap \Hcal)}$ vu dans
$\Cc^{l-s+1} \times \widetilde{G}(l-s,l-s+1)$ (on a identifié $\Hcal$
avec $\Cc^{l-s+1}$).

En résumé, quand on tranche $M_n$ avec des plans complexes $\Hcal$ de
dimension $l-s+1$, on obtient des courbes $M_n \cap
\Hcal$ dont les courbures vérifient:
$$\int \mbox{Volume}(\widetilde{(M_n \cap \Hcal)}) d \Hcal = O(\mbox{Volume}(M_n)).$$
Fixons maintenant un plan complexe $\Cc^s$ dans $\Cc ^l$ et notons
$G(l-s,l)$ l'ensemble des plans complexes de $\Cc^l$ de dimension $l-s$ qui passent
par $0$ et $\widetilde{G}(1,s)$ l'ensemble des droites du $\Cc^s$ que
l'on a fixé. L'application:
$$\Phi: \widetilde{G}(1,s) \times  G(l-s,l) \mapsto
\widetilde{G}(l-s+1,l)$$
qui envoie $(\Lcal,D)$ sur $\Lcal +D$ est définie sur un ouvert de
Zariski de $\widetilde{G}(1,s) \times  G(l-s,l)$. En
utilisant alors la formule de la coaire, on a:
$$\int \mbox{Volume}(\widetilde{(M_n \cap \Hcal(\Lcal))}) d \Lcal \mbox{ }
 dD =
O(\mbox{Volume}(M_n)),$$
où $d \Lcal$ est la mesure volume sur $\widetilde{G}(1,s)$, $dD$ celle
 de $G(l-s,l)$ et $\Hcal(\Lcal)$ est l'élément de
$\widetilde{G}(l-s+1,l)$ formé à partir de $\Lcal$ et $D$
 (i.e. $\Hcal(\Lcal)= \Lcal + D$). Par ailleurs dans l'intégrale
 ci-dessus, on ne considère que les $D$ qui font un angle au moins
 $\epsilon$ avec le $\Cc^s$ fixé (de sorte à être loin de l'ensemble
 d'indétermination de $\Phi$).

Maintenant, si on fixe $D \in G(l-s,l)$ générique au sens de la mesure, on a donc:
$$\int \mbox{Volume}(\widetilde{(M_n \cap \Hcal(\Lcal))}) d \Lcal =
O(\mbox{Volume}(M_n))$$
pour une infinité de $n$.

Dans la suite, on va utiliser ce contrôle de courbure pour construire de
bonnes îles sur les courbes $M_n \cap \Hcal(\Lcal)$ au-dessus de la
droite $\Lcal$.

\subsection{{\bf Construction des bonnes îles dans les courbes $M_n \cap \Hcal(\Lcal)$}}

Dans ce paragraphe, nous allons donner une version quantifiée de ce
que nous avons fait dans le paragraphe \ref{dim1}. Dans toute la
suite, on fixe $\Hcal(\Lcal)$ (avec $\Lcal$ générique) et on
identifiera $\Hcal(\Lcal)$ avec $\Cc^{l-s+1}$.

Soit $C_n= M_n \cap \Hcal(\Lcal)$. C'est une suite de courbes
analytiques lisses de la boule de $\Cc^{l-s+1}$ obtenue en tranchant
$B$ avec $\Hcal(\Lcal)$. Par ailleurs, on ne considèrera ici que les $\Hcal(\Lcal)$ qui
entrent dans $(1-\epsilon_0) B$. 

Reprenons ce que nous avons fait dans le paragraphe \ref{dim1}. Tout
d'abord, en utilisant le paragraphe \ref{2.1}, on a:
$$ \int_{C_n} K dV \geq - \pi \mbox{Volume}(\widetilde{C_n}).$$

Ensuite, grâce à la proposition \ref{prop2}, on sait que:
$$\mbox{Genre}(C_n \cap (1- \epsilon_0) B) \leq c(\epsilon_0) \left(
\mbox{Volume}(C_n) -  \int_{C_n} K dV \right).$$
Dans toute la suite, nous noterons $G_n$ le genre de $C_n \cap (1 -
\epsilon_0) B$ et $V_n$ le volume de $C_n $.

Maintenant, on va minorer le nombre de bonnes îles dans $C_n \cap (1 -
\epsilon_0) B$ au-dessus de $\Lcal$ en fonction de $V_n$ et du volume
de $\widetilde{C_n}$. Notons $\pi$ la projection sur $\Lcal$ en suivant la
direction $D$ (du paragraphe précédent). La projection par $\pi$ de
$C_n$ est incluse dans un carré $C$ de côté $2R$ (pour un certain $R
\in \Nn$). Si on quadrille ce carré $C \subset \Lcal$ en $4R^2k^2$
carrés égaux (de taille $\frac{1}{k}$), on a le:

\begin{lemme}

Le nombre d'îles dans $C_n \cap (1 - \epsilon_0) B$ au-dessus du quadrillage est minoré par 
$$-\frac{k}{\epsilon_k^2}(G_n+V_n) + k^2 (1- \epsilon_k) \int_{C_n \cap (1 -
  \epsilon_0 - \epsilon_k)B} \pi^{*} \omega .$$
Ici $\epsilon_k$ est une suite qui tend vers $0$ lentement et $\omega$ est la forme kählérienne
standard de $\Cc$.

\end{lemme}

\begin{proof}

Elle reprend la preuve de \cite{Det1}.

Dans un premier temps, il s'agit de modifier un peu la
courbe $C_n$ de sorte à contrôler la longueur de son bord et le nombre
de ses composantes de bord. Plus précisément, si on utilise le
paragraphe 2.1 de \cite{Det1}, on peut transformer
$C_n \cap (1- \epsilon_0)B$ en une courbe $C_n^*$ qui a son bord dans
$(1-\epsilon_0 - \epsilon_k/2)B$, qui coïncide avec $C_n$ sur
$(1-\epsilon_0 - \epsilon_k)B$ pour laquelle d'une part le nombre de composantes de bord $B_n$ est
majoré par $\frac{1}{\epsilon_k^2}(G_n + V_n)$ et d'autre part la
longueur de son bord $L_n$ est majorée par $\frac{1}{\epsilon_k} V_n$ (à des
constantes multiplicatives près que l'on oubliera).

Associé au quadrillage précédent, il y a quatre familles de $R^2k^2$ carrés deux à deux
disjoints. Dans la suite, on considère une de ces familles $Q$ et on
pave $C - Q$ en croix comme dans la figure \ref{Figure1}.

\begin{figure}[c,h]
\begin{center}

\setlength{\unitlength}{1cm}
\begin{picture}(5,3)
\includegraphics{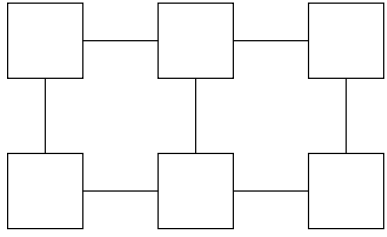}
\end{picture}

\end{center}
\caption{Pavage en croix{\label{Figure1}}. Les carrés font partie de la famille $Q$. Les croix pavent $C-Q$.}
\end{figure}

Le nombre d'îles au-dessus de $Q$ est lié à la caractéristique d'Euler
de $C_n^*-\pi^{-1}(Q)$. En effet, si $\II$ désigne l'ensemble des îles
de $\pi^{-1}(Q) \cap C_n^*$, on a $\chi(C_n^* - \pi^{-1}(Q)) \geq \chi(C_n^*) -\#\II$
(enlever une île fait chuter la caractéristique d'Euler de $1$). On
obtient donc une minoration du nombre d'îles par  $-2G_n - B_n - \chi(C_n^* -
\pi^{-1}(Q)) \geq - \frac{1}{\epsilon_k^2}(G_n+V_n)- \chi(C_n^* -
\pi^{-1}(Q))$. Il reste à majorer $\chi(C_n^* - \pi^{-1}(Q))$ pour obtenir une minoration du cardinal de $\II$.
Pour cela, on construit un graphe où chaque sommet représente une composante connexe au-dessus d'une croix, et où l'on met autant d'arêtes entre deux sommets qu'il y a d'arcs en commun dans le bord des composantes correspondantes. Dans toute la suite, on identifiera sommets et composantes connexes associées.
On obtient alors:
\begin{equation*}
\begin{split}
\chi(C_n^* - \pi^{-1}(Q))& \leq \displaystyle\sum_{ \mbox{ sommets}}{\chi(\Sigma)}-\mbox{nombre d'arêtes}\\
&\leq s-a
\end{split}
\end{equation*}

où $s$ est le nombre de sommets et $a$ le nombre d'arêtes.
La combinaison de cette relation avec la précédente, nous conduit à
une minoration du nombre d'îles au-dessus de $Q$ par $-\frac{1}{\epsilon_k^2}(G_n+V_n) + a -s$. Il nous reste donc à majorer le nombre de sommets et à minorer le
nombre d'arêtes. Pour cela, la méthode est exactement la même que dans
\cite{Det1}. Modulo un petit nettoyage des courbes
$C_n^*$, on a d'une part:
$$ s \leq R^2 S_n(C-Q)k^2(1+\epsilon_k) + \frac{k}{\epsilon_k} L_n$$
où $S_n(C-Q)$ est le recouvrement moyen de $C_n^*$ au-dessus de $C-Q$
(i.e. $S_n(C-Q)= \frac{1}{\mbox{Aire}(C-Q)} \int_{C_n^* \cap \pi^{-1}(C-Q)} \pi^{*} \omega$),
et d'autre part:
$$a \geq 2R^2 k^2 S_n(C-Q) -h k L_n,$$
où $h$ est une constante universelle (voir \cite{Det1} pour
l'obtention de ces inégalités).

On a donc trouvé une minoration du nombre d'îles au-dessus de $Q$ en 
$$-\frac{1}{\epsilon_k^2}(G_n+V_n) +R^2 k^2(1- \epsilon_k) S_n(C-Q)-
\frac{k}{\epsilon_k}L_n,$$
qui est minoré par
$$-\frac{k}{\epsilon_k^2}(G_n+V_n) +R^2 k^2(1- \epsilon_k) S_n(C-Q).$$
Maintenant, en considérant les quatre familles $Q$ de carrés dans le
quadrillage initial, on obtient une minoration du nombre d'îles par:
$$-\frac{k}{\epsilon_k^2}(G_n+V_n) + k^2(1- \epsilon_k)  \int_{C_n \cap (1 -
  \epsilon_0 - \epsilon_k)B} \pi^{*} \omega$$
qui est la minoration cherchée.

\end{proof}

On peut aussi estimer le nombre d'îles qui ne sont pas ramifiées. En
effet, grâce à un argument d'aire, on peut montrer le: 

\begin{lemme}

Le nombre d'îles non ramifiées (i.e. bonnes îles) dans $C_n \cap (1 -
\epsilon_0) B$ au-dessus du quadrillage est minoré par 
$$-\frac{k}{\epsilon_k^2}(G_n+V_n) + k^2(1- \epsilon_k)  \int_{C_n \cap (1 -
  \epsilon_0 - \epsilon_k)B} \pi^{*} \omega - k^2 \int_{C_n \cap(
  (1-\epsilon_0)B-(1-\epsilon_0-\epsilon_k)B)} \pi^{*} \omega .$$

\end{lemme}

En combinant ce lemme avec les estimées du début du paragraphe, on
obtient alors le:

\begin{lemme}{\label{lemme6}}

Le nombre de bonnes îles dans $C_n \cap (1 - \epsilon_0) B$ au-dessus du quadrillage est minoré par:
$$- \frac{k}{\epsilon_k^2} \mbox{Volume}(\widetilde{C_n}) +k^2(1- \epsilon_k) \int_{C_n \cap (1 -
  \epsilon_0 - \epsilon_k)B} \pi^{*} \omega - k^2 \int_{C_n \cap(
  (1-\epsilon_0)B-(1-\epsilon_0-\epsilon_k)B)} \pi^{*} \omega.$$

\end{lemme}

\subsection{{\bf Démonstration du critère}}

Quitte à extraire une sous-suite $\frac{[M_n]}{\mbox{Volume}(M_n)}$
converge vers un courant $T$. Le but de ce paragraphe est de démontrer
que $T$ est tissé (ou laminaire si $s=l-1$). Pour cela il suffit de le
démontrer dans $(1-\epsilon_0)B$ (avec $\epsilon_0$ petit). Dans toute
la suite, on supposera que $T_{|(1- \epsilon_0)B} \neq 0$ (sinon il
n'y a rien à faire) et que $(\pi_{D})_{*} T_{|(1 - \epsilon_0)B} \neq
0$ (où $\pi_D$ est la projection sur le $\Cc^s$ que l'on avait fixé en
suivant la direction générique $D$ du paragraphe \ref{3.1}).

Rappelons que pour démontrer le critère, nous allons utiliser le théorème de N. Sibony et
P. M. Wong (voir \cite{SW}). Plus précisément, l'énoncé que nous
allons utiliser est le (voir \cite{DS} ou lemme 3.7 de \cite{Di}): 

\begin{theoreme*}  \cite{SW}

Soit $\Dcal$ une famille de droites dans $\Cc^{s}$ passant par un
point $a$ et soit $B'$ la boule unité de $\Cc^s$ centrée en $a$.

Supposons que $\Hcal_{2s-2}(\Dcal) \geq \frac{1}{2} $ où
$\Hcal_{2s-2}$ est la mesure volume sur l'ensemble des droites de
$\Cc^s$ qui passent par $a$. Soit $g$ une
fonction qui est holomorphe au voisinage de $a$ ainsi que sur les droites de $\Dcal$
intersectées avec $rB'$. Alors $g$ se prolonge en une fonction holomorphe dans la boule $c r B'$
(ici $c$ est indépendante de $g$ et $\Dcal$). De plus, on a:
$$\sup_{b \in c r B'} |g(b) - g(a)| \leq \sup_{b \in \Dcal \cap r B'}
|g(b) - g(a)|.$$

\end{theoreme*}

Reprenons les notations des paragraphes précédents. Si on considère un
  point $x$ de $\Cc^s$ et une droite $L$ de $G(1,s)$ (qui est
  l'ensemble des droites de $\Cc^s$ qui passent par $0$), nous noterons
  toujours $\Hcal(x+L)$ le plan complexe de dimension $l-s+1$ associé
  (i.e. $\Hcal(x+L)=x + L +D$). 

Notons maintenant $m(x,L)$ le nombre de bonnes îles dans $M_n \cap \Hcal(x+L)
\cap (1-\epsilon_0)B$ au-dessus du carré $S$ où $S \subset (x+L)$ est le carré de
taille $\frac{1}{k}$ centré en $x$. 

Voici le plan de la démonstration. Dans un premier paragraphe, on va
minorer $\int m(x,L) dL \mbox{ } dx$. Dans le second, on verra que la
combinaison de cette minoration avec le théorème de N. Sibony et
P. M. Wong permettra de construire beaucoup de bonnes îles (de
dimension $s$) dans les sous-ensembles analytiques $M_n$. Enfin dans
le dernier, on démontrera le caractère tissé de $T$.  

\subsubsection{{\bf Minoration de $\int m(x,L) dL \mbox{ } dx$}}    

Considérons l'application $\Phi: \Cc^s \times G(1,s) \mapsto
\widetilde{G}(1,s)$ qui envoie $(x,L)$ sur $x+L$. En utilisant la
formule de coaire, on a:
$$\int \mbox{apJac}_{4(s-1)} \Phi m(x,L) dL \mbox{ } dx=\int
\int_{\Lcal} m(y,L(\Lcal)) dy \mbox{ } d \Lcal,$$
où $dy$ est la mesure de Lebesgue sur $\Cc$ et $L(\Lcal)$ est la direction
de $\Lcal$.

Le jacobien $\mbox{apJac}_{4(s-1)} \Phi$ ne jouera aucun rôle dans la
suite: on pourra donc l'oublier. Nous sommes donc ramenés à minorer:
$$\int \int_{\Lcal} m(y,L(\Lcal)) dy \mbox{ } d \Lcal.$$
Pour faire cette minoration, nous allons utiliser les estimées du paragraphe précédent.

Fixons une droite $\Lcal$. L'image de la boule unité $B$
par $\pi_D$ (où $\pi_D$ est le projection sur $\Cc^s$ en suivant la
direction $D$) dans $\Lcal$ est incluse dans un carré de longueur
$2R$ (pour un certain $R \in \Nn$). Soit
maintenant $\Qcal_k$ le quadrillage de ce carré en $4R^2 k^2$ carrés égaux
de taille $\frac{1}{k}$ et $dy_0$ la mesure de Lebesgue sur un de ces carrés $C_0$. Si $y_0$ est dans $C_0$, on notera
$y_i(y_0)$ les $4R^2 k^2$ points des carrés du quadrillage qui dans
ceux-ci ont la même
position que $y_0$ dans $C_0$. Ainsi, on a:
$$\int_{\Lcal}  m(y,L(\Lcal)) dy = \int_{C_0} \sum_{i=1}^{4R^2 k^2} m(y_i(y_0), L(\Lcal)) dy_0.$$
Autrement dit, par le lemme \ref{lemme6},

$$\int m(x,L) dL \mbox{ } dx = \int \int_{C_0} \sum_{i=1}^{4R^2k^2}
m(y_i(y_0), L(\Lcal)) dy_0 \mbox{ } d \Lcal$$

$$\geq (1- \epsilon_k) \frac{k^2}{k^2} \int \left( \int_{M_n
  \cap \Hcal(\Lcal) \cap (1-\epsilon_0- \epsilon_k)B} \pi^{*}
\omega \right) d \Lcal - \frac{a_{n,k}}{k^2}$$
avec
$$a_{n,k}= \frac{k}{ \epsilon_k^2} \int  \mbox{Volume}(\widetilde{M_n \cap
  \Hcal(\Lcal)}) d \Lcal+ b_{n,k},$$
où
$$b_{n,k}= k^2 \int \int_{M_n \cap \Hcal(\Lcal) \cap
  ((1-\epsilon_0)B-(1-\epsilon_0-\epsilon_k)B)} \pi^{*} \omega \mbox{
  }  d \Lcal.$$ 
Dans ces expressions, $\pi$ est la projection sur $\Lcal$ (dans $\Hcal(\Lcal)$) induite par $D$ et $\omega$
est la forme volume sur $\Lcal$.

Maintenant, grâce au paragraphe \ref{3.1}, on a:
$$\int \mbox{Volume}(\widetilde{M_n \cap \Hcal(\Lcal) }) d \Lcal=O(\mbox{Volume}(M_n)).$$
Autrement dit,
$$\int m(x,L) dL \mbox{ } dx \geq (1 -\epsilon_k) \int \mbox{Volume}(\pi
(M_n \cap \Hcal(\Lcal) \cap (1-\epsilon_0- \epsilon_k)B)) d \Lcal
-\frac{1}{k \epsilon_k^2} O(\mbox{Volume}(M_n)) - \frac{b_{n,k}}{k^2}.$$
En particulier, si on note $n(y,\Lcal)$ le nombre
d'antécédents de $y \in \Lcal$ par $\pi$ dans $M_n \cap \Hcal(\Lcal)
\cap (1-\epsilon_0- \epsilon_k)B$, on a:
$$\int m(x,L) dL \mbox{ } dx \geq (1 -\epsilon_k) \int \int_{\Lcal} n(y,\Lcal) dy \mbox{ } d\Lcal -\frac{1}{k \epsilon_k^2} O(\mbox{Volume}(M_n))- \frac{b_{n,k}}{k^2},$$
d'où (toujours en oubliant le jacobien de $\Phi$),
$$\int m(x,L) dL \mbox{ } dx \geq (1 -\epsilon_k) \int n(x) dL \mbox{
} dx -\frac{1}{k \epsilon_k^2} O(\mbox{Volume}(M_n))- \frac{b_{n,k}}{k^2},$$
où $n(x)$ est le nombre de relevés du point $x$ dans $M_n \cap
(1-\epsilon_0 - \epsilon_k)B$. Mais si $\epsilon_0$ est suffisamment
générique, alors l'inégalité précédente reste vraie si $n(x)$ est le
nombre de relevés du point $x$ dans $M_n \cap (1-\epsilon_0 )B$. Dans
toute la suite $n(x)$ désignera cette dernière quantité. De plus,
toujours parce que $\epsilon_0$ est suffisamment générique, on a:
$$\frac{b_{n,k}}{k^2} \leq \epsilon_k \int n(x) dL \mbox{ } dx.$$

En particulier, comme $m(x,L) \leq n(x)$ on en déduit que:
$$0 \leq \int(n(x)-m(x,L))dL \mbox{ } dx \leq \epsilon_k
\mbox{Volume}(M_n),$$
avec $\epsilon_k$ qui tend vers $0$ lentement.

\subsubsection{{\bf Construction des bonnes îles}}

Grâce à l'estimée obtenue dans le paragraphe précédent, on va pouvoir utiliser le théorème de N. Sibony
et P. M. Wong. Pour cela soit ($n$ est fixé):
$$X_k= \{ x, \mbox{ } n(x)(1- \epsilon_k') \geq m(x,L) \mbox{ pour un ensemble } \Sigma_x
\mbox{ de } L
\mbox{ de mesure supérieure à } \epsilon_k' \}.$$
Il s'agit de montrer que cet ensemble
est petit.

Si on reprend l'estimée précédente, on a:
$$\epsilon_k \mbox{Volume}(M_n) \geq \int_{X_k} \int_{\Sigma_x} (n(x)-m(x,L))dL
\mbox{ } dx \geq \epsilon_k' \int_{X_k} \int_{\Sigma_x} n(x) dL \mbox{ } dx$$
Ce qui implique que l'on a une majoration de la forme:
$$\int_{X_k} n(x) dx \leq \epsilon_k  \mbox{Volume}(M_n),$$
pourvu que la suite $(\epsilon_k')^2$ tende vers $0$ moins vite que
$\epsilon_k$.

Considérons un point $x$ hors de $X_k$ qui n'est pas dans les
valeurs critiques de $\pi_{D |M_n}$. L'ensemble des droites $L$
pour lesquelles $n(x)(1-\epsilon_k') \leq m(x,L)$ est de mesure au
moins $1- \epsilon_k'$. Notons $x_1, \cdots , x_{n(x)}$
les relevés du point $x$ dans $M_n \cap (1 - \epsilon_0)B$. On dira que $L$ se relève bien en
$x_i$ s'il existe une bonne île dans $M_n \cap (1 - \epsilon_0)B$ qui contient
$x_i$ et qui est au-dessus du carré de taille
$\frac{1}{k}$ centré en $x$ dans $x+L$. Comme $x$
n'est pas dans $X_k$, le nombre
de $x_i$ pour lesquels la mesure de $L$ qui se relèvent bien en $x_i$ est
majoré par $\frac{1}{2}$ est au plus égal à $\epsilon_k' n(x)$ (à des
constantes multiplicatives près). En utilisant maintenant le théorème
de N. Sibony et P. M. Wong, on en déduit que l'on peut construire $(1-\epsilon_k')n(x)$
relevés du cube de taille $\frac{c}{k}$ centré en $x$ dans $M_n \cap
(1 - \epsilon_0)B$. La restriction de $\pi_D$ à ces relevés est un
biholomorphisme sur le cube de taille $\frac{c}{k}$ centré en $x$.

\subsubsection{{\bf Caractère tissé de la limite}}

Dans ce paragraphe nous allons montrer que la limite $T$ est
tissée dans $(1- \epsilon_0)B$.

Rappelons que $(\pi_{D})_{*} T_{|(1 - \epsilon_0)B} \neq
0$ (car la direction $D$ est choisie générique). L'image de la boule $B$ par la projection $\pi_D$ sur $\Cc^s$ est
incluse dans le cube $C$ centré en $0$ de taille $cR$ (pour un certain $R \in
\Nn$). Considérons le découpage de $C$ en cubes de taille
$\frac{c}{2k}$. Nous appellerons bonnes îles les préimages $P$ par $\pi_D$
de ces cubes dans $M_n \cap (1 - \epsilon_0)B$ pour lesquelles la restriction de $\pi_D$ à
$P$ est un biholomorphisme sur le cube du découpage qui lui
correspond. D'après l'estimée sur le volume de $X_k$ et le nombre de
branches inverses au-dessus d'un point hors de $X_k$, on peut minorer
le nombre de bonnes îles au-dessus du quadrillage par $\left( \frac{2k}{c}
\right)^{2s}(1- \epsilon_k')(\mbox{Volume}(\pi_D(M_n \cap (1- \epsilon_0)B)) - \epsilon_k \mbox{Volume}(M_n))$. Maintenant comme $(\pi_{D})_{*} T_{|(1
  - \epsilon_0)B} \neq 0$, on en déduit que le nombre de bonnes îles
est minoré par $\left( \frac{2k}{c} \right)^{2s}(1 - \epsilon_k)
\mbox{Volume}(\pi_D(M_n \cap (1 - \epsilon_0)B))$.

Soit $T_{k,n}$ le courant défini par
$T_{k,n}=\frac{1}{ \mbox{Volume}(M_n)}
\displaystyle\sum_{ \mbox{ bonnes îles}}[\Gamma]$. Le courant
$T_{k,n}$ peut aussi s'écrire $T_{k,n}=\int [\Gamma] d \nu_{k,n}(\Gamma)$ où
$\nu_{k,n}$ est une mesure sur l'espace métrique compact des graphes de
dimension $s$ au-dessus des cubes du quadrillage.

Si on note $T_n=\frac{[M_n \cap (1-\epsilon_0) B]}{\mbox{Volume}(M_n)}$ alors on a:
$$\int T_{k,n} \wedge \pi_D^{*} (\omega^{s}) \geq
(1-\epsilon_k)\int{T_{n} \wedge \pi_D^{*}( \omega^{s})},$$
d'où,
$$\int{(T_{n}-T_{k,n}) \wedge \pi_D^{*} (\omega^{s})} \leq \epsilon_k.$$
La suite de mesures $\nu_{k,n}$ converge vers une mesure $\nu_k$
(quitte à extraire une sous-suite) ce qui implique que $T_{k,n}$
converge vers $T_k= \int [\Gamma] d \nu_k(\Gamma) $ qui est donc uniformément
tissé au-dessus de chaque cube du quadrillage (et uniformément
laminaire si $s=l-1$). Par ailleurs, on a toujours l'estimée:
$$\int{(T_{| (1-\epsilon_0) B}-T_{k}) \wedge \pi_D^{*} (\omega^{s})} \leq \epsilon_k,$$
avec $T_{| (1-\epsilon_0) B}-T_k \geq 0$ par construction.
Si on raffine de plus en plus le quadrillage (i.e. si $k$ augmente),
$T_k$ croît vers un courant $T_{\infty}$ qui est tissé (ou laminaire
si $s=l-1$). De plus $T_{\infty} \leq T_{| (1-\epsilon_0) B}$ et
$\int{(T_{| (1-\epsilon_0) B}-T_{\infty}) \wedge \pi_D^{*} (\omega^{s})} \leq
0$.

Maintenant, si on prend une autre direction $D'$ générique par
rapport à $T_{\infty}$ et telle que $(\pi_{D'})_{*}(T_{| (1-\epsilon_0) B}) \neq 0$ (où
$\pi_{D'}$ désigne la projection associée à $D'$), on construit de
même un courant $T_{\infty}' \leq T_{| (1-\epsilon_0) B}$ qui est supérieur à
$T_{\infty}$ et qui vérifie  $\int{(T_{| (1-\epsilon_0) B} -T_{\infty}') \wedge
  \pi_{D'}^{*} (\omega^s)} = 0$. En itérant ce procédé on finit par avoir $T_{|
  (1- \epsilon_0) B}=T_{\infty}'$, c'est-à-dire que $T$ est tissé
dans $ (1- \epsilon_0) B$ (ou laminaire si $s=l-1$). C'est ce que l'on
voulait démontrer.

$ $

Henry de Thélin

Université Paris-Sud (Paris 11)

Mathématique, Bât. 425

91405 Orsay

France


\begin{thebibliography}{00}



\bibitem{BLS} E. Bedford, M. Lyubich et J. Smillie, \textit{Polynomial
  diffeomorphisms of $\Cc^2$. IV: The measure of maximal entropy and
  laminar currents}, Invent. Math., {\bf 112} (1993), 77-125.
\bibitem{Bi} E. Bishop,  \textit{Conditions for the analyticity of certain sets}, Michigan Math. J., {\bf 11} (1964), 289-304.
\bibitem{C} S. Cantat, \textit{Dynamique des automorphismes des
  surfaces K3}, Acta Math., {\bf 187} (2001), 1-57. 
\bibitem{Ch} E. M. Chirka, \textit{Complex analytic sets}, Kluwer
  Academic Publishers, Dordrecht (1989). 
\bibitem{Det1} H. de Thélin, \textit{Sur la laminarité de certains
courants}, Ann. Sci. Ecole Norm. Sup., {\bf 37} (2004), 304-311.
\bibitem{Det2} H. de Thélin, \textit{Un phénomène de concentration de
  genre}, Math. Ann., {\bf 332} (2005), 483-498.
\bibitem{Di} T.-C. Dinh, \textit{Suites d'applications méromorphes
  multivaluées et courants laminaires}, J. Geom. Anal., {\bf 15} (2005), 207-227.
\bibitem{DS} T.-C. Dinh et N. Sibony, \textit{Dynamique des
  applications d'allure polynomiale}, J. Math. Pures Appl., {\bf 82}
  (2003), 367-423.
\bibitem{Du1} R. Dujardin, \textit{Laminar currents in $\Pp^2$},
  Math. Ann., {\bf 325} (2003), 745-765.
\bibitem{DuvSi} J. Duval et N. Sibony,  \textit{Polynomial convexity, rational convexity, and currents}, Duke Math. J., {\bf 79} (1995), 487-513.
\bibitem{Fe} H. Federer, \textit{Geometric measure theory}, Springer
  Verlag (1969). 
\bibitem{LS} R. Langevin et T. Shifrin, \textit{Polar varieties and
  integral geometry}, Amer. J. Math., {\bf 104} (1982), 553-605.
\bibitem{L} P. Lelong, \textit{Propriétés métriques des variétés analytiques complexes définies par une équation}, Ann. Sci. Ecole Norm. Sup., {\bf 67} (1950), 393-419.
\bibitem{SW} N. Sibony et P. M. Wong, \textit{Some results on global
  analytic sets}, Séminaire Pierre Lelong-Henri Skoda, Lecture Notes
  in Math., {\bf 822} (1980), 221-237.

\end{thebibliography}
\end{document}